\documentclass[12pt,a4paper]{article}
\usepackage[width=18cm,height=23cm]{geometry}
\usepackage{amsmath,amssymb,amsfonts,amsthm,enumerate}
\providecommand{\keywords}[1]{\textbf{\textit{Keywords: }}#1}
\newcommand{\cond}[1]{\quad{\scriptstyle{(#1)}}}
\newcommand{\bc}{\mathbb{C}}
\newcommand{\bq}{\mathbb{Q}}

\newcommand{\bz}{\mathbb{Z}}
\newtheorem{Thm}{Theorem}[section]%
\newtheorem{Lem}[Thm]{Lemma}%

\newtheorem{Cor}[Thm]{Corollary}
\newtheorem{Prop}[Thm]{Proposition}%
\theoremstyle{definition}
\newtheorem{Def}[Thm]{Definition}%
\newtheorem{Exa}[Thm]{Example}%
\newtheorem{Quest}[Thm]{Question}%
\usepackage{indent first}

\title{\bf On nonsingularity of circulant matrices}
\author{Zhangchi Chen}
\date{\today}

\begin{document}
\maketitle

\abstract{In Communication theory and Coding, it is expected that certain circulant matrices having $k$ ones and $k+1$ zeros in the first row are nonsingular. We prove that such matrices are always nonsingular when $2k+1$ is either a power of a prime, or a product of two distinct primes. For any other integer $2k+1$ we construct circulant matrices having determinant $0$. The smallest singular matrix appears when $2k+1=45$. The possibility for such matrices to be singular is rather low, smaller than $10^{-4}$ in this case.}



\section{Introduction}
\label{sect-intro}

We begin with definitions and some classical results on circulant matrices, taken from (\cite{Ingleton-1956} Section 1).
\begin{Def} A \emph{circulant matrix} $C(a_0,\dots,a_{n-1})$ is an $n\times n$ matrix of the form
\[
\begin{bmatrix}
    a_0       & a_1  & \dots & a_{n-1} \\
    a_{n-1} & a_0  & \dots & a_{n-2} \\
    \vdots & \vdots & \ddots & \vdots \\
    a_1 & a_2 & \dots & a_0
\end{bmatrix},
\]
where $a_0,\dots,a_{n-1}\in\bq$. It is \emph{unital} if $a_0,\dots,a_{n-1}\in\{0,1\}$.
\end{Def}

The determinant of a circulant matrix $C(a_0,\dots,a_{n-1})$ can be calculated by
\[
\det C(a_0,\dots,a_{n-1})=\prod\limits_{j=0}^{n-1}\big(a_0+a_1\,\omega_j+a_2\,\omega_j^2+\dots a_{n-1}\,\omega_j^{n-1}\big),
\]
in terms of the $n$-th root of unity
\[
\omega_j:=e^{\frac{2\pi ij}{n}},\cond{i^2=-1}.
\]

A proof can be found in (\cite{Golub-Van Loan-2013} Theorem 4.8.2).

\begin{Def} The polynomial
\[
f(x)=a_0+a_1\,x+a_2\,x^2+\dots+a_{n-1}\,x^{n-1}\in\bq[x]
\]
is called the \emph{associated polynomial} of $C(a_0,\dots,a_{n-1})$. It is called \emph{unital} if all its coefficients $a_0,\dots,a_{n-1}$ are in $\{0,1\}$.
\end{Def}

Note that the relation between circulant matrices and their associated polynomials is not `1-1'. For example the identity matrix of any size has $f(x)=1$. In fact, a circulant matrix is determined by both its associated polynomial and its size.
 
\begin{Prop}\label{prop-common-roots} A circulant matrix $C(a_0,\dots,a_{n-1})$ is nonsingular if and only if its associated polynomial $f(x)$ and $x^n-1$ share no common roots.
\end{Prop}
\proof By what precedes:
\[
\det C(a_0,\dots,a_{n-1})=\prod\limits_{j=0}^{n-1}f(\omega_j).\qedhere
\]
\endproof

\begin{Def} A circulant matrix $C(a_0,\dots,a_{n-1})$ is called \emph{$r$-recurrent} with $r$ a proper divisor of $n$, if $a_j=a_{j'}$ whenever $j\equiv j'$ mod $r$. It is called \emph{non-recurrent} if it is not $r$-recurrent for any proper divisor $r$ of $n$.
\end{Def}

If $C(a_0,\dots,a_{n-1})$ is $r$-recurrent for some proper divisor $r$ of $n$, then it has the same rank as $C(a_0,\dots,a_{r-1})$ because
\[
C(a_0,\dots,a_{n-1})
=
\begin{bmatrix}
    C(a_0,\dots,a_{r-1})       & C(a_0,\dots,a_{r-1}) & \dots & C(a_0,\dots,a_{r-1}) \\
    C(a_0,\dots,a_{r-1})       & C(a_0,\dots,a_{r-1}) & \dots & C(a_0,\dots,a_{r-1}) \\
    \vdots & \vdots & \ddots & \vdots \\
    C(a_0,\dots,a_{r-1})      & C(a_0,\dots,a_{r-1}) & \dots & C(a_0,\dots,a_{r-1})
\end{bmatrix}.
\]
Hence an invertible circulant matrix must be non-recurrent. This paper mainly studies circulant unital matrices having $k$ ones and $k+1$ zeros in the first row. They are always non-recurrent.

One motivation to study such matrices comes from Communication and Coding where $2k+1$ input signals are `mixed' by these matrices \cite[Theorem~13, Remark~5]{Wan-2017}, \cite[Theorem~29]{Wan-2018}. Experts hope such processes to be invertible. Another motivation comes from \cite{Newton-1954} on summability of polydiagonal matrices for periodic sequences of zeros and ones. The following two questions are equivalent:

\begin{Quest}\label{quest-sing} For fixed $k\in\bz_{\geqslant 1}$, does there exist a circulant unital matrix having $k$ ones and $k+1$ zeros in its first row which is singular?
\end{Quest}

\begin{Quest} For fixed $k\in\bz_{\geqslant 1}$, does there exist a unital polynomial $f(x)\in\bz[x]$ such that $\deg f(x)\leqslant 2k$, $f(1)=k$ and $f(x)$ shares a common root with $x^n-1$?
\end{Quest}

In this paper we propose a complete answer.

\begin{Thm}\label{thm:pe} If $2k+1=p^e$ for some prime $p$, then such matrices are always nonsingular.
\end{Thm}

\begin{Thm}\label{thm:pq} If $2k+1=p\,q$ for two distinct primes $3\leqslant p<q$, then such matrices are always nonsingular.
\end{Thm}

\begin{Thm}\label{thm:pqr} If $2k+1=p\,q\,r$ where $3\leqslant p<q$ are two distinct primes and $r\geqslant 3$ is an odd integer, then there exist some singular matrices of such a type.
\end{Thm}

\begin{Cor}\label{thm:45} The first singular example appears when $(p,q,r)=(3,5,3)$, i.e. $2k+1=45$. For example if
\begin{align*}
E_{22}:= \{ &0, 9, 18, 27, 36,\\
& 3, 12, 21, 30, 39,\\
&1, 16, 31,\\
&2, 17, 32,\\
&4, 19, 34,\\
&5, 20, 35\}.
\end{align*}
and if $a_j=1$ for $j\in E_{22}$ and $a_j=0$ for $j\notin E_{22}$, then $C(a_0,\dots,a_{44})$ is singular.
\end{Cor}
\proof Let $f_{E_{22}}(x)\in\bz[x]$ be the associated polynomial of this matrix. Then
\begin{align*}
f_{E_{22}}(x) =\sum\limits_{j\in E_{22}}x^j&=(1+x^3)\,(1+x^9+x^{18}+x^{27}+x^{36})+(x+x^2+x^4+x^5)\,(1+x^{15}+x^{30})\\
&=(1+x^3)\,{\textstyle{\frac{x^{45}-1}{x^9-1}}}+(x+x^2+x^4+x^5)\,\textstyle{\frac{x^{45}-1}{x^{15}-1}}
\end{align*}
annihilates $e^{\frac{2\pi i}{45}}$, a root of $x^{45}-1$.\qed

\section{Cyclotomic polynomials}
\label{sect-basic}

By Proposition \ref{prop-common-roots}, $C(a_0,\dots,a_{n-1})$ is singular if and only if $f(x)$ and $x^n-1$ have a common root. Thus we should study the irreducible and the unital factors of $x^n-1$, namely the cyclotomic polynomials and the fundamental recurrent polynomials (Definition \ref{def-frec}).

\begin{Def} For any $n\in\bz_{\geqslant 1}$, the \emph{$n$-th cyclotomic polynomial} is defined as
\[
\Phi_n(x):=\prod\limits_{
1\leqslant d\leqslant n\atop
\gcd(d,n)=1}
\big(x-e^{\frac{2\pi i d}{n}}\big)\in\bc[x].
\]
\end{Def}

We recall the following properties (\cite[VI. 3]{Lang-2002}).
\begin{itemize}
\item (Degree) $\deg \Phi_n(x)=\varphi(n)$ where $\varphi$ is Euler's totient function.
\item (Integer coefficients) $\Phi_n(x)\in\bz[x]$.
\item (Irreducibility) $\Phi_n(x)$ is irreducible over $\bz[x]$.
\item (Factorization) $x^n-1=\prod\limits_{r|n}\Phi_r(x)$.
\item For any two distinct $n,m\in\bz_{\geqslant 1}$, $\Phi_n(x)$ and $\Phi_m(x)$ do not divide one another.
\item For a prime number $p$ and for $e\in\bz_{\geqslant 1}$, we have $\Phi_{p^e}(x)=1+x^{p^{e-1}}+x^{2p^{e-1}}+\dots+x^{p^e-p^{e-1}}$. Hence $\Phi_{p^e}(1)=p$. Note that $x^{p^e}-1=\Phi_1(x)\,\Phi_p(x)\dots\Phi_{p^e}(x)$.
\end{itemize}

\begin{Thm} A circulant matrix of size $n\times n$ is singular if and only if its associated polynomial $f(x)$ is divisible by $\Phi_d(x)$ for some $d$ dividing $n$.\qed
\end{Thm}

\section{Case where $2k+1$ is a power of a prime}
\label{sect-casepe}

From (\cite[IV.~2, Corollary~2.2]{Lang-2002}), recall the Gauss Lemma: if a monic polynomial $f(x)\in\bz[x]$ is equal to $h(x)\,g(x)$ with $g(x)\in\bz[x]$ monic and $h(x)\in\bq[x]$, then $h(x)\in\bz[x]$ has integer coefficients.


\begin{Lem}\label{lem:p} Let $f(x)$ be a unital polynomial, and let $p$ be a prime number not dividing $f(1)$. Then $\Phi_{p^e}(x)$ does not divide $f(x)$ for any $e\in\bz_{\geqslant 0}$.
\end{Lem}

\proof When $e=0$, since $f(1)\in\bz_{\geqslant 1}$ is nonzero, $\Phi_1(x)=x-1$ cannot divide $f(x)$. When $e\in\bz_{\geqslant 1}$, suppose $g(x)=\Phi_{p^e}(x)$ divides $f(x)$ for some $e\in\bz_{\geqslant 1}$, so that by the Gauss Lemma, $f(x)=h(x)\,g(x)$ with $h(x)\in\bz[x]$. Hence $f(1)=h(1)\,g(1)=h(1)\,p$, contradicting the hypothesis that $p$ does not divide $f(1)$.\qed

\proof[Proof of Theorem~\ref{thm:pe}] Write $2k+1=p^e$ with $p$ an odd prime and $e\in\bz_{\geqslant 1}$. Let $f(x)\in\bz[x]$ be a unital polynomial with $\deg f(x)<p^e-1$ and $f(1)=k$. Remember $x^{p^e}-1=\Phi_1(x)\,\Phi_p(x)\dots\Phi_{p^e}(x)$. Since $p$ does not divide $k$, by Lemma \ref{lem:p}, $\Phi_{p^{e'}}(x)$ does not divide $f(x)$ for any $e'\in\bz_{\geqslant 0}$. Thus $f(x)$ shares no common root with $x^{p^e}-1$, i.e. the corresponding circulant matrix is nonsingular.\endproof

\section{Recurrent decompositions}

Next we introduce the unital factors of $x^n-1$: the fundamental recurrent polynomials. The following notation is taken from Ingleton \cite{Ingleton-1956}.

\begin{Def}\label{def-frec} For any proper divisor $r$ of $n$, the \emph{fundamental $r$-recurrent polynomial with respect to $n$} is
\[
G(n,r;x):=\textstyle{\frac{x^n-1}{x^r-1}}=1+x^r+\dots+x^{n-r}\in\bz[x].
\]
\end{Def}

The name comes from the following fact: for any $r$-recurrent $n\times n$ circulant matrix, the associated polynomial $f(x)$ is a multiple of $G(n,r;x)$.

We have
\begin{equation}\label{decompo-G}
G(n,r;x)=\prod\limits_{d|n}\Phi_d(x)\Big/\prod\limits_{d|r}\Phi_d(x)=\prod\limits_{d|n\atop d\nmid r}\Phi_d(x),
\end{equation}
and in particular, $\Phi_n(x)$ divides $G(n,r;x)$.

For example $G(45,9;x)=\Phi_{45}(x)\,\Phi_{15}(x)\,\Phi_5(x)$.

In the previous section, note that we have decomposed 
\begin{align*}
f_{E_{22}}(x)
=(1+x^3)\,G(45,9;x)+(x+x^2+x^4+x^5)\,G(45,15;x).
\end{align*}
In general, if $f(x)\in\bq[x]$ with $\deg f(x)<n$ can be decomposed as $\sum\limits_{r|n\atop r<n}h_r(x)\,G(n,r;x)$ for some $h_r(x)\in\bc[x]$, then clearly $\Phi_n(x)$ divides $f(x)$.

In this section we want to establish the converse. For any $f(x)\in\bq[x]$ with $\deg f(x)<n$, which is divisible by $\Phi_n(x)$, we want to decompose it as $\sum\limits_{r|n\atop r<n}h_r(x)\,G(n,r;x)$.

\begin{Def} \em{(Recurrent decomposition)} Fix $n\in\bz_{\geqslant 1}$ and let $f(x)\in\bq[x]$ with $\deg f(x)<n$. Let $p_1<\dots<p_m$ be all the distinct prime factors of $n$. If
\begin{align*}
f(x)=\sum\limits_{j=1}^m h_{\frac{n}{p_j}}(x)\,G(n,{\textstyle{\frac{n}{p_j}}};x),
\end{align*}
for some $h_{\frac{n}{p_j}}(x)\in\bq[x]$ with $\deg h_{\frac{n}{p_j}}(x)<\frac{n}{p_j}$, then we call $\big(h_{\frac{n}{p_1}}(x),\dots,h_{\frac{n}{p_m}}(x)\big)$ a $(p_1,\dots,p_m)$-recurrent decomposition of $f(x)$ with respect to $n$. We call such a decomposition unital if each $h_{\frac{n}{p_j}}(x)$ is unital.

\end{Def}

For example $(x+x^2+x^4+x^5,1+x^3)$ is a unital $(3,5)$-recurrent decomposition of $f_{E_{22}}(x)$ with respect to $45$.

The following Theorem is equivalent to \cite[Proposition~3.1]{Ingleton-1956}. 

\begin{Thm}{\em(Existence of a recurrent decomposition)}\label{thm-rd} Let $f(x)\in \bq[x]$ with $\deg f(x)<n$. Suppose $\Phi_n(x)$ divides $f(x)$. Let $p_1<\dots<p_m$ be all the distinct prime factors of $n$. Then $f(x)$ admits a $(p_1,\dots,p_m)$-recurrent decomposition with respect to $n$.
\end{Thm}

\proof Using (\ref{decompo-G}), we have
\begin{align*}
\gcd\Big(\{G(n,{\textstyle{\frac{n}{p_j}}};x)\colon\,j=1,\dots,m\}\Big)
=
\gcd\Big(\{\prod\limits_{d|n,\,d\nmid\frac{n}{p_j}}\Phi_d(x)\colon\,j=1,\dots,m\}\Big)=\Phi_n(x).
\end{align*}
Applying the Euclidean Algorithm we find some $g_{\frac{n}{p_j}}(x)\in\bq[x]$ such that
\[
\Phi_n(x)=\sum\limits_{j=1}^mg_{\frac{n}{p_j}}(x)\,G(n,\textstyle\frac{n}{p_j};x).
\]
Suppose $f(x)=q(x)\,\Phi_n(x)$ where $q(x)\in\bq[x]$. Then $f(x)=\sum\limits_{j=1}^mq(x)\,g_{\frac{n}{p_j}}(x)\,G(n,\frac{n}{p_j};x)$. We now need to change the multiplicators $q(x)\,g_{\frac{n}{p_j}}(x)$ to meet the degree bounds.

The quotient
\[
\frac{x^{n}-1}{G(n,\frac{n}{p_j};x)}=x^{\frac{n}{p_j}}-1
\]
is a polynomial of degree $\frac{n}{p_j}$. By Euclidean division, there exists a unique $h_{\frac{n}{p_j}}(x)\in\bq[x]$ with $\deg h_{\frac{n}{p_j}}(x)<\frac{n}{p_j}$ such that
$q(x)\,g_{\frac{n}{p_j}}(x)\equiv h_{\frac{n}{p_j}}(x)$ mod $x^{\frac{n}{p_j}}-1$. Then
\[
f(x)=\sum\limits_{j=1}^m q(x)\,g_{\frac{n}{p_j}}(x)\,G(n,\textstyle{\frac{n}{p_j}};x)\equiv\sum\limits_{j=1}^mh_{\frac{n}{p_j}}(x)\,G(n,\frac{n}{p_j};x) ~\text{mod}~ x^n-1.
\] Since $\deg f(x)<n$ we conclude that $f(x)=\sum\limits_{j=1}^mh_{\frac{n}{p_j}}(x)\,G(n,\frac{n}{p_j};x)$.\qed

When $f(x)$ is unital, we expect the existence of a unital decomposition. This is not true in general. A counterexample where $m=3$ will be constructed at the end of the section. The case where $m=2$ is proved by Ingleton \cite[4.1]{Ingleton-1956}. We prove a stronger version here.

We write $\bz_{[0,d]}[x]$ for the set of polynomials having coefficients in $\{0,1,\dots,d\}$.
\begin{Thm}\label{thm:decompo-unital} {\em(Recurrent decomposition in $\bz_{[0,d]}[x]$)} Let $p<q$ be two distinct primes. Let $n=p^{e_1}\,q^{e_2}$ with $e_1,\,e_2\in\bz_{\geqslant 1}$. Let $f(x)\in\bz_{[0,d]}[x]$ with $\deg f(x)<n$ be a polynomial which is divisible by $\Phi_n(x)$. Then $f(x)$ admits a $(p,q)$-recurrent decomposition in $\bz_{[0,d]}[x]$ with respect to $n$. 

In particular, when $d=1$, a unital decomposition exists. 
\end{Thm}
\proof By Theorem \ref{thm-rd}, there exist $h_{\frac{n}{p}}(x),\,h_{\frac{n}{q}}(x)\in\bq[x]$ with $\deg h_{\frac{n}{p}}(x)<\frac{n}{p},\,\deg h_{\frac{n}{q}}(x)<\frac{n}{q}$ such that $f(x)=h_{\frac{n}{p}}(x)\,G(n,\frac{n}{p};x)+h_{\frac{n}{q}}(x)\,G(n,\frac{n}{q};x)$. Write
\begin{align*}
f(x) &=\sum\limits_{j=0}^{n-1}a_j\,x^j =\sum\limits_{s=0}^{\frac{n}{p\,q}-1}x^s\Big(\sum\limits_{l=0}^{p\,q-1}a_{\frac{l\,n}{p\,q}+s}\,x^{\frac{l\,n}{p\,q}}\Big),\\
h_{\frac{n}{p}}(x)&=\sum\limits_{u=0}^{\frac{n}{p}-1}b_u\,x^u=\sum\limits_{s=0}^{\frac{n}{p\,q}-1}x^s\Big(\sum\limits_{l=0}^{q-1}b_{\frac{l\,n}{p\,q}+s}\,x^{\frac{l\,n}{p\,q}}\Big),\\
h_{\frac{n}{q}}(x)&=\sum\limits_{v=0}^{\frac{n}{q}-1}c_v\,x^v=\sum\limits_{s=0}^{\frac{n}{p\,q}-1}x^s\Big(\sum\limits_{l=0}^{p-1}c_{\frac{l\,n}{p\,q}+s}\,x^{\frac{l\,n}{p\,q}}\Big).
\end{align*}
For every $s=0,\dots,\frac{n}{p\,q}-1$, we group the coefficients as
\begin{align}
A_s&:=\{a_{\frac{l\,n}{p\,q}+s}:l=0,\dots,p\,q-1\}\subset\{0,\dots,d\},\nonumber\\
B_s&:=\{b_{\frac{l\,n}{p\,q}+s}:l=0,\dots,q-1\}\subset\bq,\label{Bs}\\
C_s&:=\{c_{\frac{l\,n}{p\,q}+s}:l=0,\dots,p-1\}\subset\bq.\nonumber
\end{align}
Define $e_s:=\min B_s\in\bq$ and introduce
\begin{align*}
g(x)&:=\sum\limits_{s=0}^{\frac{n}{p\,q}-1}e_s\,x^s.
\end{align*}
Using
\begin{align*}
G(\textstyle\frac{n}{p},\frac{n}{p\,q};x)=\sum\limits_{l=0}^{q-1}x^{\frac{l\,n}{p\,q}}, && G(\textstyle\frac{n}{q},\frac{n}{p\,q};x)=\sum\limits_{l=0}^{p-1}x^{\frac{l\,n}{p\,q}},
\end{align*}
we may subtract and modify the multiplicators as
\begin{align*}
h'_{\frac{n}{p}}(x)&:=h_{\frac{n}{p}}(x)-g(x)\,G(\textstyle\frac{n}{p},\frac{n}{p\,q};x)=\sum\limits_{s=0}^{\frac{n}{p\,q}-1}x^s\Big(\sum\limits_{l=0}^{q-1}(b_{\frac{l\,n}{p\,q}+s}-e_s)\,x^{\frac{l\,n}{p\,q}}\Big),\\
h'_{\frac{n}{q}}(x)&:=h_{\frac{n}{q}}(x)+g(x)\,G(\textstyle\frac{n}{q},\frac{n}{p\,q};x)=\sum\limits_{s=0}^{\frac{n}{p\,q}-1}x^s\Big(\sum\limits_{l=0}^{p-1}(c_{\frac{l\,n}{p\,q}+s}+e_s)\,x^{\frac{l\,n}{p\,q}}\Big).
\end{align*}
Since $G(\frac{n}{p},\frac{n}{p\,q};x)\,G(n,\frac{n}{p};x)=G(\frac{n}{q},\frac{n}{p\,q};x)\,G(n,\frac{n}{q};x)$, we see that $\big(h'_{\frac{n}{p}}(x),h'_{\frac{n}{q}}(x)\big)$ is another decomposition of $f(x)$.
So we can suppose from the beginning that $\min B_s=0$ for all $s$. 

Given two sets of rational numbers $X,\,Y\subset\bq$, define their sum and difference by
\[
X+ Y:=\{x+ y:x\in X,\,y\in Y\},\qquad X-Y:=\{x-y:x\in X,y\in Y\}.
\]
We have $X\subset (X+Y)-Y$. From
\begin{align*}
h_{\frac{n}{p}}(x)\,G(n,\textstyle\frac{n}{p};x)&=\sum\limits_{s=0}^{\frac{n}{p\,q}-1}x^s\Big(\sum\limits_{l=0}^{q-1}b_{\frac{l\,n}{p\,q}+s}\,\big(x^{\frac{l\,n}{p\,q}}+x^{\frac{(l+q)\,n}{p\,q}}+\dots+x^{\frac{(l+(p-1)q)\,n}{p\,q}}\big)\Big),\\
h_{\frac{n}{q}}(x)\,G(n,\textstyle\frac{n}{q};x)&=\sum\limits_{s=0}^{\frac{n}{p\,q}-1}x^s\Big(\sum\limits_{l=0}^{p-1}c_{\frac{l\,n}{p\,q}+s}\,\big(x^{\frac{l\,n}{p\,q}}+x^{\frac{(l+p)\,n}{p\,q}}+\dots+x^{\frac{(l+(q-1)p)\,n}{p\,q}}\big)\Big),\\
\end{align*}
we deduce that each coefficient $a_{\frac{l_1\,n}{p\,q}+s}\in A_s$, $l_1\in\{0,\dots,p\,q-1\}$ is a sum of the form
\[
a_{\frac{l_1\,n}{p\,q}+s}=b_{\frac{l_2\,n}{p\,q}+s}+c_{\frac{l_3\,n}{p\,q}+s}, 
\]
where $l_2\in\{0,\dots,q-1\}$, $l_3\in\{0,\dots,p-1\}$ such that $l_1\equiv l_2$ mod $q$ and $l_1\equiv l_3$ mod $p$. Thus by the Chinese Reminder Theorem, 
\[
B_s+C_s=A_s\subset\{0,1,\dots,d\}.
\]
Since $\min B_s=0$ we have $C_s=\{0\}+C_s\subset B_s+C_s\subset\{0,1,\dots,d\}$ and $B_s\subset (A_s-C_s)\cap \bq_{\geqslant 0}=\{-d,-d+1,\dots,d\}\cap \bq_{\geqslant 0}=\{0,1,\dots,d\}$, i.e. both $h_{\frac{n}{p}}(x)$ and $h_{\frac{n}{q}}(x)$ are in $\bz_{[0,d]}[x]$.
\qed

\begin{Exa}
If $n$ has more than two distinct prime factors, Theorem \ref{thm:decompo-unital} may be false, i.e. a unital polynomial $f(x)$ divisible by $\Phi_n(x)$ may not admit any unital recurrent decomposition. A counterexample appears when $n=105=3\times 5\times 7$ and 
\begin{align*}
f(x)&=x^5+x^6+x^{10}+x^{25}+x^{27}+x^{35}+x^{40}+x^{48}+x^{50}+x^{65}+x^{69}+x^{70}+x^{80}+x^{85}+x^{95}+x^{100}\\
&=(1+x^{5}+x^{10}+x^{15}+x^{25}+x^{30})\,G(105,35;x)+x^6\,G(105,21;x)-G(105,15;x).
\end{align*}
Here $f(x)$ is unital and $f(1)=16$. However if we suppose $f(x)=h_{35}(x)\,G(105,35;x)+h_{21}(x)\,G(105,21;x)+h_{15}(x)\,G(105,15;x)$ for some unital $h_{35}(x),\,h_{21}(x)$ and $h_{15}(x)$, then $16=3\,h_{35}(1)+5\,h_{21}(1)+7\,h_{15}(1)$, where $h_{35}(1),\,h_{21}(1),\,h_{15}(1)\in\bz_{\geqslant0}$. The only solutions $\big(h_{35}(1),h_{21}(1),h_{15}(1)\big)$ to this equation are $(0,2,2)$ and $(1,0,3)$, hence either $h_{21}(1)$ or $h_{15}(1)$ is 0, i.e. either $h_{21}$ or $h_{15}$ is 0 since they are unital. In the first case $e^{\frac{2\pi i}{15}}$ is a root of $h_{35}(x)\,G(105,35;x)+h_{21}(x)\,G(105,21;x)$ but not a root of $f(x)$. In the second case $e^{\frac{2\pi i}{21}}$ is a root of $h_{35}(x)\,G(105,35;x)+h_{15}(x)\,G(105,15;x)$ but not a root of $f(x)$. These contradictions prove that $f(x)$ admits no unital recurrent decomposition.

\begin{table}[!htb]
\caption{The values of each polynomial at 3 special points}
\begin{center}
\begin{tabular}{c|c|c|c|c}\hline
& $G(105,35;x)$ & $G(105,21;x)$ & $G(105,15;x)$ & $f(x)$\\\hline
$x=e^{\frac{2\pi i}{35}}$ & 3&0 &0 &$-3e^{\frac{8\pi i}{7}}$\\\hline
$x=e^{\frac{2\pi i}{21}}$ & 0& 5&0 & $5e^{\frac{6\pi i}{7}}$\\\hline
$x=e^{\frac{2\pi i}{15}}$ & 0& 0& 7&-7 \\\hline
\end{tabular}
\end{center}
\label{default}
\end{table}%
\end{Exa}
\section{Case where $2k+1$ is a product of two distinct primes}
\label{sect-casepq}
\proof[Proof of Theorem~\ref{thm:pq}] Write $2k+1=p\,q$, where $p,\,q$ are two distinct primes. Suppose $C(a_0,\dots,a_{2k})$ is singular. Let $f(x)\in\bz_{[0,1]}[x]$ be its corresponding unital polynomial. Since neither $p$ nor $q$ divides $k$, by Lemma \ref{lem:p} we know none among $\Phi_1(x),\,\Phi_p(x)$ and $\Phi_q(x)$ divides $f(x)$. So $\Phi_{p\,q}(x)$ divides $f(x)$. By Theorem \ref{thm:decompo-unital} we have
\[
f(x)=h_q(x)\,G(p\,q,q;x)+h_p(x)\,G(p\,q,p;x),
\]
for some $h_q(x),\,h_p(x)$ unital. As before, we write
\begin{align*}
f(x)=\sum\limits_{j=0}^{p\,q-1}a_j\,x^j, && h_q(x)=\sum\limits_{u=0}^{q-1}b_u\,x^u, && h_p(x)=\sum\limits_{v=0}^{p-1}c_v\,x^v,
\end{align*}
and write
\begin{align*}
A_0&:=\{a_l:l=0,\dots,p\,q-1\},\\
B_0&:=\{b_l:l=0,\dots,q-1\},\\
C_0&:=\{c_l:l=0,\dots,p-1\}.
\end{align*}
We have $k=f(1)=h_q(1)\,q+h_p(1)\,p$. Since neither $p$ nor $q$ divides $k$, we know $h_q(1)\neq0$, $h_p(1)\neq0$, hence $1\in B_0$, $1\in C_0$. By the Chinese Reminder Theorem we know $B_0+C_0=A_0$, hence $2\in A_0$, contradicting the hypothesis that $f(x)$ is unital.\qed



\section{Other cases: constructing a singular matrix}
\label{sect-casepqr}

\proof[Proof of Theorem~\ref{thm:pqr}] In all other cases we can write $2k+1=p\,q\,r$ where $p,\,q$ are two distinct primes and $r\geqslant3$ is an odd integer. We may assume $p<q$ and $p\leqslant r$ by choosing $p,\,q$ as the first and the second smallest prime factors of $2k+1$. To construct a singular circulant matrix of our type, it suffices to find a unital polynomial $f(x)\in\bz[x]$ with $\deg f(x)<p\,q\,r,\,f(1)=k$ such that $\Phi_{p\,q\,r}(x)$ divides $f(x)$. Since $k=\frac{p\,q\,r-1}{2}\geqslant\frac{3p\,q-1}{2}>p\,q$, we have unique $a,\,b\in\bz_{\geqslant 1}$ with $b\leqslant p-1$ such that $a\,p+b\,q=k$. Define

\begin{align*}
f(x):&=(1+x^{q\,r}+x^{2q\,r}+\dots+x^{b\,q\,r-q\,r})\,G(p\,q\,r,p\,r)(x)+\sum\limits_{j\in R_a}x^j\,G(p\,q\,r,q\,r)(x),
\end{align*}
where $R_a\subset \{0,1,\dots,q\,r-1\}\backslash\{0,r,\dots,q\,r-r\}$ is a set of $a$ elements. Such $R_a$ exists if and only if $a\leqslant q\,r-q$. This is true since

\begin{align*}
q\,r-q-a &= q\,r-q-\textstyle{\frac{k-b\,q}{p}}\\
 & = q\,r-q-\textstyle{\frac{p\,q\,r-1-b\,q}{2p}}\\
 &=\textstyle{\frac{q\,r}{2}}-q+\frac{1+b\,q}{2p}\\
  &\geqslant\textstyle{\frac{r-2}{2}}\,q>0.
\end{align*}

Clearly we have $f(x)\in\bz[x]$ and $\deg f(x)<p\,q\,r$. Since $\Phi_{p\,q\,r}(x)$ divides both $G(p\,q\,r,p\,r;x)$ and $G(p\,q\,r,q\,r;x)$ it also divides $f(x)$. Moreover $f(1)=b\,G(p\,q\,r,p\,r)(1)+a\,G(p\,q\,r,q\,r)(1)=b\,q+a\,p=k$. The condition $R_a\cap\{0,r,\dots,q\,r-r\}=\emptyset$ ensures that $f(x)$ is unital.\qed

The smallest $k$ in this case is $22$, when $2k+1=45=3^2\times 5$. Here we take $p=r=3,\,q=5$. Note that $22=4\times 3+2\times 5$. In this case we can construct a singular matrix as above, by taking $R_4=\{1,2,4,5\}$. Then $f(x)=(1+x^3)\,(1+x^9+x^{18}+x^{27}+x^{36})+(x+x^2+x^4+x^5)\,(1+x^{15}+x^{30})$ represents a $45\times 45$ singular circulant unital matrix with 22 ones in its first row. This reveals how we constructed $E_{22}$ in Corollary \ref{thm:45}.

\section{The number of singular matrices when $k=22$}

Finally we count the number of such singular matrices and estimate the probability of a circulant matrix of our type to be singular. We study the easiest case when $k=22$. There are in total $\binom{45}{22}=4116715363800$ choices of such matrices. Suppose $f(x)$ is a corresponding polynomial of a singular circulant matrix of our type. Then there exists some $r$ dividing $45$ such that $\Phi_r(x)$ divides $f(x)$. Since neither $3$ nor $5$ divides $f(1)$, by Lemma \ref{lem:p} we know none among $\Phi_1(x),\,\Phi_3(x),\,\Phi_5(x)$ and $\Phi_9(x)$ divides $f(x)$. Thus either $\Phi_{15}(x)$ or $\Phi_{45}(x)$ divides $f(x)$.

\subsection{Case (1): $\Phi_{45}(x)$ divides $f(x)$}

Theorem \ref{thm:decompo-unital} guarantees the existence of unital recurrent decompositions of $f(x)$, but not the uniqueness. However, we can always control the ambiguity.

\begin{Thm}\label{thm:ambiguity} Suppose $n=p^{e_1}\,q^{e_2}$ has two distinct prime factors $p,\,q$. If $f(x)\in\bq[x]$ with $\deg f(x)$ divisible by $\Phi_n(x)$, then for any two $(p,q)$-recurrent decompositions $\big(h_{\frac{n}{p}}(x),h_{\frac{n}{q}}(x)\big)$ and $\big(h'_{\frac{n}{p}}(x),h'_{\frac{n}{q}}(x)\big)$ of $f(x)$ with respect to $n$ such that
\begin{align*}
\deg h_{\frac{n}{p}}(x)<\frac{n}{p},&&\deg h'_{\frac{n}{p}}(x)<\frac{n}{p},&&\deg h_{\frac{n}{q}}(x)<\frac{n}{q},&&\deg h'_{\frac{n}{q}}(x)<\frac{n}{q},
\end{align*}
there exists some $\delta(x)\in\bq[x]$ with $\deg \delta(x)<\frac{n}{p\,q}$ such that
\begin{align}\label{eqn:ambiguity}
h'_{\frac{n}{p}}(x)=h_{\frac{n}{p}}(x)-\delta(x)\,G(\textstyle{\frac{n}{p}},\frac{n}{p\,q};x), && h'_{\frac{n}{q}}(x)=h_{\frac{n}{q}}(x)+\delta(x)\,G(\textstyle{\frac{n}{q}},\frac{n}{p\,q};x).
\end{align}
\end{Thm}

\proof By the definition of $(p,q)$-recurrent decompositions with respect to $n$, we have
\begin{align}\label{eqn-diff}
\big(h_{\frac{n}{p}}(x)-h'_{\frac{n}{p}}(x)\big)\,G(n,\textstyle{\frac{n}{p}};x)=\big(h'_{\frac{n}{q}}(x)-h_{\frac{n}{q}}(x)\big)\,G(n,\textstyle{\frac{n}{q}};x).
\end{align}
Note that
\begin{align*}
G(n,\textstyle{\frac{n}{p}};x)&=\Phi_n(x)\prod\limits_{j=1}^{e_2}\Phi_{n\,q^{-j}}(x)=\Phi_n(x)\,G(\textstyle{\frac{n}{q}},\frac{n}{p\,q};x),\\
G(n,\textstyle{\frac{n}{q}};x)&=\Phi_n(x)\prod\limits_{j=1}^{e_1}\Phi_{n\,p^{-j}}(x)=\Phi_n(x)\,G(\textstyle{\frac{n}{p}},\frac{n}{p\,q};x).
\end{align*}
Hence $G(\frac{n}{q},\frac{n}{p\,q};x)$ and $G(\frac{n}{p},\frac{n}{p\,q};x)$ share no nontrivial common factors. Dividing both sides of (\ref{eqn-diff}) by $\Phi_n(x)$, we get
\begin{align*}
\big(h_{\frac{n}{p}}(x)-h'_{\frac{n}{p}}(x)\big)\,G(\textstyle{\frac{n}{q}},\frac{n}{p\,q};x)=\big(h'_{\frac{n}{q}}(x)-h_{\frac{n}{q}}(x)\big)\,G(\textstyle{\frac{n}{p}},\frac{n}{p\,q};x).
\end{align*}
So there exists some $\delta(x)\in\bq[x]$ such that 
\begin{align*}
h_{\frac{n}{p}}(x)-h'_{\frac{n}{p}}(x)=\delta(x)\,G(\textstyle{\frac{n}{p}},\frac{n}{p\,q};x)&&h'_{\frac{n}{q}}(x)-h_{\frac{n}{q}}(x)=\delta(x)\,G(\textstyle{\frac{n}{q}},\frac{n}{p\,q};x).
\end{align*}
We have $\deg \delta(x)\leqslant \max\{h_{\frac{n}{p}}(x),h'_{\frac{n}{p}}(x)\}-\deg G(\frac{n}{p},\frac{n}{p\,q};x)< \frac{n}{p}-(\frac{n}{p}-\frac{n}{p\,q})=\frac{n}{p\,q}$.\qed

\begin{Def}\label{def:ambiguity} As in the proof of Theorem \ref{thm:decompo-unital}, for any $(p,q)$-recurrent decomposition $\big(h_{\frac{n}{p}}(x),h_{\frac{n}{q}}(x)\big)$ of $f(x)$ with respect to $n$, write
\[
h_{\frac{n}{p}}(x)=\sum\limits_{u=0}^{\frac{n}{p}-1}b_u\,x^u,
\]
and for any $s=0,\dots,\frac{n}{p\,q}-1$ define
\[
B_s:=\{b_{\frac{l\,n}{p\,q}+s}\colon l=0,\dots,q-1\}.
\]
We call a $(p,q)$-recurrent decomposition {\em $p$-uniformized} if $\min B_s=0$ for any $s=0,\dots,\frac{n}{p\,q}-1$.
\end{Def}

\begin{Cor}\label{cor:p-uniformized} Suppose $n=p^{e_1}\,q^{e_2}$ has two distinct prime factors $p,\,q$. Let $f(x)\in\bq[x]$ with $\deg f(x)<n$ be divisible by $\Phi_n(x)$. Then among all the $(p,q)$-recurrent decompositions of $f(x)$ with respect to $n$, there exists a unique $p$-uniformized one. Moreover, if $f(x)$ is unital, its $p$-uniformized decomposition is also unital.
\end{Cor}
\proof (Existence) For any $(p,q)$-recurrent decomposition $\big(h_{\frac{n}{p}}(x),h_{\frac{n}{q}}(x)\big)$ of $f(x)$ with respect to $n$, define $e_s:=\min B_s\in\bq$, $g(x):=\sum\limits_{s=0}^{\frac{n}{p\,q}-1}e_s\,x^s$. Then $\big(h_{\frac{n}{p}}(x)-g(x)\,G(\frac{n}{p},\frac{n}{p\,q};x),h_{\frac{n}{q}}(x)+g(x)\,G(\frac{n}{q},\frac{n}{p\,q};x)\big)$ is $p$-uniformized.

(Uniqueness) If there are two $p$-uniformized $(p,q)$-recurrent decompositions $\big(h_{\frac{n}{p}}(x),h_{\frac{n}{q}}(x)\big)$ and $\big(h'_{\frac{n}{p}}(x),h'_{\frac{n}{q}}(x)\big)$ of $f(x)$ with respect to $n$, by Theorem \ref{thm:ambiguity} there exists some $\delta(x)\in\bq[x]$ with $\deg \delta(x)<\frac{n}{p\,q}$ satisfying (\ref{eqn:ambiguity}).

Write $B_s,\,B'_s$ as a collection of coefficients of $h'_{\frac{n}{p}}(x)$ and $h'_{\frac{n}{p}}(x)$ respectively, as in Definition~\ref{def:ambiguity}. Write $\delta(x)=\sum\limits_{s=0}^{\frac{n}{p\,q}-1}\delta_s\,x^s$. From (\ref{eqn:ambiguity}) we have $\min B'_s=\min B_s-\delta_s$. So $\delta_s=0$ for all $s=0,\dots,\frac{n}{p\,q}-1$. We conclude that these two decompositions are equivalent.

(Unital) In Theorem \ref{thm:decompo-unital} we have already constructed a $p$-uniformized unital $(p,q)$-recurrent decomposition when $f(x)$ is unital.
\endproof

By Corollary \ref{cor:p-uniformized} we have a unique $3$-uniformized unital $(3,5)$-recurrent decomposition of $f(x)$ with respect to $45$, i.e.
\[
f(x)=h_{15}(x)\,G(45,15;x)+h_9(x)\,G(45,9;x),
\]
where $h_{15}(x),\,h_9(x)\in\bz[x]$ are unital and if we write
\begin{align*}
f(x)=\sum\limits_{j=0}^{44}a_j\,x^j, && h_{15}(x)=\sum\limits_{u=0}^{14}b_u\,x^u, && h_9(x)=\sum\limits_{v=0}^{8}c_v\,x^v,
\end{align*}
for $s=0,1,2$ we may define 
\begin{align*}
A_s&:=\{a_{3l+s}:l=0,\dots,14\}\subset\{0,1\},\\
B_s&:=\{b_{3l+s}:l=0,\dots,4\}\subset\{0,1\},\\
C_s&:=\{c_{3l+s}:l=0,\dots,2\}\subset\{0,1\}.
\end{align*}
Here we have $\min B_s=0$ for any $s=0,1,2$. Moreover $22=f(1)=3h_{15}(1)+5h_9(1)$ where $h_{15}(1),\,h_9(1)\in\bz_{\geqslant 0}$, hence $h_{15}(1)=4,\,h_9(1)=2$.

Each $f(x)$ corresponds to a unique choice of $b_0,\dots,b_{14},\,c_0,\dots,c_8\in\{0,1\}$ satisfying the following conditions:
\begin{itemize}
\item $\sum\limits_{u=0}^{14} b_u=4,\,\sum\limits_{v=0}^{8}c_v=2$,
\item $\min B_s=0$,
\item if $1\in B_s$ then $C_s=\{0\}$.
\end{itemize}

Case (1.1): there is only 1 $s$ such that $1\in C_s$. There are $\binom{3}{1}\binom{3}{2}\binom{10}{4}=1890$ choices.

Case (1.2): there are 2 $s$ such that $1\in C_s$. There are $\binom{3}{2}\binom{3}{1}^2\binom{5}{4}=135$ choices.

\smallskip
Hence there are 2025 choices of unital $f(x)$ divisible by $\Phi_{45}(x)$, corresponding to 2025 singular matrices.

\subsection{Case (2): $\Phi_{15}(x)$ divides $f(x)$}
Write $f^{(15)}(x):=\sum\limits_{j=0}^{14}d_jx^j\in\bz[x]$ where $d_j:=a_j+a_{j+15}+a_{j+30}\in\{0,1,2,3\}$. In fact it is the residue of $f(x)$ divided by $x^{15}-1$, hence divisible by $\Phi_{15}(x)$. Such $f^{(15)}(x)$ corresponds to $\prod\limits_{j=0}^{14}\binom{3}{d_j}$ unital $f(x)$. For any $f^{(15)}(x)$ divisible by $\Phi_{15}(x)$, by Corollary \ref{cor:p-uniformized} it admits a unique $3$-uniformized unital $(3,5)$-recurrent decomposition with respect to $15$:
\[
f^{(15)}(x)=h_5(x)\,G(15,5;x)+h_3(x)\,G(15,3;x),
\]
where $h_5(x),\,h_3(x)\in\bq[x]$ and if we write
\begin{align*}
h_5(x)=\sum\limits_{u=0}^{4}b_u\,x^u, && h_3(x)=\sum\limits_{v=0}^{2}c_v\,x^v,
\end{align*}
define 
\begin{align*}
A_0&:=\{d_{l}:l=0,\dots,14\}\subset\{0,1,2,3\},\\
B_0&:=\{b_{l}:l=0,\dots,4\},\\
C_0&:=\{c_{l}:l=0,\dots,2\}.
\end{align*}
Then $\min B_0=0$. By the Chinese Reminder Theorem we have $B_0+C_0=A_0$. Thus $C_0=\{0\}+C_0\subset B_0+C_0\subset\{0,1,2,3\}$ and $B_0\subset [0,+\infty)\cap (A_0+\{0,-1,-2,-3\})=\{0,1,2,3\}$. Moreover we have $f^{(15)}(1)=22=3h_5(1)+5h_3(1)$ where $h_3(1), h_5(1)\in\bz_{\geqslant 0}$, hence $h_5(1)=4,\,h_3(1)=2$.

Each $f^{(15)}(x)$ corresponds to a unique choice of $b_0,\dots,b_4,\,c_0,\dots,c_2\in\{0,1,2,3\}$ satisfying the following conditions:
\begin{enumerate}
\item $\sum\limits_{u=0}^{4} b_u=4,\,\sum\limits_{v=0}^{2}c_v=2$,
\item $\min B_0=0$,
\item $\max B_0+\max C_0\leqslant 3$.
\end{enumerate}

We can write down a complete list of possible values of $(b_u),\,(c_v)$ up to permutations. Note that
\begin{itemize}
\item each $f^{(15)}(x)$ is uniquely determined by values of $(b_u),\,(c_v)$;
\item for each $f^{(15)}(x)=\sum\limits_{j=0}^{14}d_j\,x^j$ there are $\prod\limits_{j=0}^{14}\binom{3}{d_j}=3^{\#\{j\colon d_j=1,2\}}$ choices of unital $f(x)$.
\end{itemize} 

\begin{table}[!htb]
\caption{Possible values of $(b_0,\dots,b_4),\,(c_0,c_1,c_2)$ up to permutations}
\begin{center}
\begin{tabular}{c|c|c}\hline
Type of values & Permutations & Choices of $f(x)$ for each $f^{(15)}(x)$\\\hline
$(0,0,0,2,2),\,(0,1,1)$ & $\binom{5}{2}\binom{3}{2}=60$ & $3^6\times 3^2=3^8$ \\\hline
$(0,0,1,1,2),\,(0,1,1)$ & $\binom{5}{1}\binom{4}{2}\binom{3}{2}=90$ & $3^6\times 3^5=3^{11}$ \\\hline
$(0,1,1,1,1),\,(0,0,2)$ & $\binom{5}{4}\binom{3}{1}=15$ & $3^8\times 3^1=3^9$ \\\hline
$(0,1,1,1,1),\,(0,1,1)$ & $\binom{5}{4}\binom{3}{2}=15$ & $3^6\times 3^8=3^{14}$ \\\hline
\end{tabular}
\end{center}
\label{default}
\end{table}%

By summing up $(\text{Permutations})\times\big(\text{Choices of $f(x)$ for each $f^{(15)}(x)$}\big)$ we get $88376670$ choices of unital $f(x)$ divisible by $\Phi_{15}(x)$, corresponding to $88376670$ singular matrices.

\subsection{Double counts: both $\Phi_{15}(x)$ and $\Phi_{45}(x)$ divide $f(x)$}
Suppose a unital $f(x)$ with $\deg f(x)<45$ is divisible by both $\Phi_{15}(x)$ and $\Phi_{45}(x)$. As in case (1) there exist some unital $h_{15}(x),\,h_9(x)$ such that
\[
f(x)=h_{15}(x)\,G(45,15;x)+h_9(x)\,G(45,9;x),
\]
with $h_{15}(1)=4$. Take value at $\zeta_{15}:=e^{\frac{2\pi i}{15}}$. Since $f(x)$ is divisible by $\Phi_{15}(x)$ we have $f(\zeta_{15})=0$. We also have $G(45,15;\zeta_{15})=3,\,G(45,9;\zeta_{15})=0$. Hence $h_{15}(\zeta_{15})=0$. We conclude that $h_{15}(x)$ is divisible by $\Phi_{15}(x)$. It is also unital and $\deg h_{15}(x)<15$. By Theorem~\ref{thm:decompo-unital}, $h_{15}(x)$ admits a unital $(3,5)$-recurrent decomposition with respect to $15$:
\[
h_{15}(x)=h_5(x)\,G(15,5;x)+h_3(x)\,G(15,3;x),
\]
for some unital $h_5(x),\,h_3(x)\in\bz[x]$ with $\deg h_5(x)<5$, $\deg h_3(x)<3$.

We have $4=h_{15}(1)=3h_5(1)+5h_3(1)$, where $h_5(1),\,h_3(1)\in\bz_{\geqslant 0}$. However there is no solution $\big(h_5(1),h_3(1)\big)$ to this equation. We conclude that there is no unital $f(x)$ with $\deg f(x)<45,\,f(1)=22$ and $f(x)$ divisible by both $\Phi_{15}(x)$ and $\Phi_{45}(x)$.

We have $88378695$ singular unital circulant matrices having exactly $22$ ones in their first rows. The possibility of a unital circulant matrix in our type being singular is about $2.15\times10^{-5}<10^{-4}$. This algorithm can be generalized to all $n$ with only 2 distinct prime factors.

\end{document}